\documentclass{article}
\usepackage{latexsym}

\newcommand{\bc}[2]{{{#1}\choose{#2}}}

\begin{document}

\author{B\'{e}la Bajnok, 
Steven B. Damelin, Jenny Li, and 
Gary L. Mullen}
\title{A Constructive Finite Field Method for Scattering Points on the
Surface of $d$-Dimensional Spheres}
\date{April 27, 2001}
\maketitle

\begin{abstract} 
We use solutions to quadratic forms in $d+1$ variables over finite fields
 to scatter points on the surface of the unit sphere $S^d$, 
$d\geq 1$. 
Applications are given for spherical $t$ designs and generalized $s$ energies.
\end{abstract}

\section{Introduction} 
\setcounter{equation}{0}

In this exploratory article, we present
a constructive method for scattering points on the surface of
$d$ dimensional spheres which we believe is new and of interest. 
Indeed, the problem of uniformly distributing
points on spheres is an interesting and difficult problem 
with vast applications in fields as diverse as crystallography,
approximation theory,
computational complexity, molecular structure, and electrostatics.
We refer the interested reader to \cite{ABG}, \cite{BBP}, \cite{CS},
\cite{DT}, \cite{EH}, \cite{GraTic}, \cite{KM}, \cite{KS1},
\cite{RSZ}, \cite{RSZ1}
and the
many references cited therein
for a comprehensive account of this fascinating subject.

Let $d \geq 1$ be an integer and $S^d$ the $d$-dimensional unit sphere 
given by the set of all real solutions to the equation
\begin{equation}\label{eq 1.1}
x_1^2 + \cdots + x_{d+1}^2=1. 
\end{equation}
In two dimensions $(d=1)$ the problem is easily reduced to uniformly
distributing $N$ points on a circle and the vertices of the regular $N$-gon
provide an obvious answer. For $d\geq 2$ the problem becomes much more
difficult; in fact, there are numerous criteria for uniformity, resulting
in different   optimal
configurations, see \cite{HS}. 
 
In this article, we will describe a constructive method to
scatter points on the surface of $S^d$ using finite fields. After
describing our construction, we will apply it to two measures 
of uniformity on $S^d$, namely spherical $t$-designs and generalized $s$
energy. 

\newpage

{\bf Finite Field Construction}
\medskip

For an odd prime $p$, let $F_p$ denote the
finite field of integers modulo $p$. 
Consider the quadratic form given by (1.1) over $F_p$.
The number $N=N(d,p)$ of solutions  of this form is well known and is given by

\begin{equation}\label{eq 1.2}
N(d,p)=
\left\{ \begin{array}{cc}
p^d -p^{(d-1)/2} \eta ((-1)^{(d+1)/2}) & {\rm if}\, d \mbox{  is odd}  \\
p^d + p^{d/2} \eta((-1)^{d/2}) & {\rm if}\, d \mbox{  is even}
\end{array} \right. 
\end{equation}
where $\eta$ is the quadratic character defined on $F_p$ by $\eta(0)=0,$
$\eta(a )=1 $ if $a$ is a square in 
$F_p$, and $\eta(a)=-1$ if $a$ is a non-square in $F_p$, see 
\cite[Theorems 6.26 and 6.27]{LN}. 
 
\medskip

Given a solution vector 
\[
X= (x_1, \dots, x_{d+1}),\, x_i\in F_p,\, 1\leq i\leq d+1,
\] 
we may assume without loss of generality that the points $x_i$ are scaled
so that 
they  
are centered around the origin and are contained in the set
\[
\{-(p-1)/2,...,(p-1)/2\}.
\]
More precisely, if $x_i\in X$, define
\begin{displaymath}
w_i=
\left\{ \begin{array}{cc}
x_i & {\rm if}\, x_i\in \{0,...,(p-1)/2\}  \\
x_i-p & {\rm if}\, x_i\in \{(p+1)/2,...,p-1\}.
\end{array} \right.
\end{displaymath}
Then $w_i\in \{-(p-1)/2,...,(p-1)/2\}$ and 
the scaled vector
\[
W= (w_1, \dots, w_{d+1}),\,1\leq i\leq d+1
\]
solves (1.1) if and only if $X$ solves (1.1). 
 
Denoting by $||;||$ the usual Euclidean metric, we
multiply each solution vector $W$ by $\frac{1}{||W||}$.
Clearly each of these normalized points 
is now on the surface of the unit sphere $S^d$. 
Use of the finite field $F_p$ for larger primes $p$ 
provides a method to increase the number $N$ of points that are placed on the 
surface of $S^d$ for any fixed $d\geq 1$. For increasing values of 
$p$, we obtain an
increasing number $N = O(p^d)$ of points scattered on the surface of the
unit sphere 
$S^d$; in particular, as $p \rightarrow \infty$ through all odd primes, it
is clear that 
$N \rightarrow \infty$.
For each prime $p$ and integer $d \geq 1$,
we will henceforth denote the set of points arising from our
finite field construction by $X=X(d,p)$.
\medskip

Let us now describe the point set $X$ produced by the finite field
construction and provide some explicit examples for small values of $p$ and
$d$.  In each case, we may start with a well chosen set $V=V(d,p)$ of vectors. Then, in order to construct the full set of points $X(d,p)$,
we need to
consider all points obtained from $V$ by taking $\pm 1$ times the entry in
each
coordinate, and by permuting the coordinates of each vector, in all
possible ways.  For small values of $d$ and $p$, this construction is
summarized in the following table.

$$\begin{array}{||c|c||c|c||} \hline \hline
d & p & N(d,p) & V(d,p)  \\ \hline \hline
1 & 3 & 4 & \{(1,0)\}  \\ \hline
1 & 5 & 4 & \{(1,0)\}  \\ \hline
1 & 7 & 8 & \{(1,0), \frac{1}{\sqrt{2}} (1,1) \}  \\ \hline \hline
2 & 3 & 6 & \{(1,0,0)\}  \\ \hline
2 & 5 & 30 & \{(1,0,0),\frac{1}{\sqrt{2}} (2,1,1)\}  \\ \hline
2 & 7 & 42 & \{(1,0,0),\frac{1}{\sqrt{2}} (1,1,0),\frac{1}{\sqrt{22}}
(3,3,2)\}  \\ \hline \hline
\end{array}$$

\medskip

Observe that for $p=3,5,7$ and $d=1$, our construction gives the optimal
solution, namely the 
the vertices of the regular $N$-gon. This, however, is not the case for
$p>7$. 
\medskip

The remainder of this paper is organized as follows.  In Section 2, we will
discuss spherical $t$-designs. More precisely, we will prove that for
any $d\geq 2$ and any odd prime $p$, our set of points forms a spherical
$3$-design. Moreover, given any odd positive integer $k$, our set of points
forms a spherical design of index $k$.
Section 3 introduces the notion of generalized $s$ energies.
Here we numerically compare
our results with well known theoretical and numerical bounds
for these latter quantities using results of
Wagner, Kuijlaars and Saff, and Rakhmanov, Saff, and Zhou, see  \cite{W1}, \cite{W2},
\cite{KS1}, \cite{RSZ}, and \cite{RSZ1}. Finally, 
in Section 4, we briefly discuss a natural extension of our method to
finite fields of prime power orders.       
Appendices A and B contain numerical data to
illustrate the effectiveness of our construction. The computer programs
used to make these calculations can be obtained from the  authors. Throughout, $C$ will denote a positive constant which may take on different values at different times. 
\bigskip

\section{Spherical $t$-designs}

\setcounter{equation}{0}

In this section, we will study how well our points are distributed
using spherical $t$-designs, a notion first introduced by Delsarte,
Goethals, and Seidel
in \cite{DGS}.
 
\bigskip

{\bf Definition 2.1}\,  
A finite set $X$ of points on the $d$-sphere $S^d$ is a \emph{spherical
$t$-design} or a \emph{spherical design of strength $t$}, if for every
polynomial $f$ of total degree $t$ or less, the average value of $f$ over
the whole sphere is equal to the arithmetic average of its values on $X$.
If this only holds for 
homogeneous polynomials of degree $t$, then $X$ is called a \emph{spherical
design of index $t$}.
\bigskip

In other words, $X$ is of index $t$ if the Chebyshev-type quadrature formula
\begin{equation} \label{eq:quad}
\frac{1}{\sigma_d(S^d)} \int_{S^d}f({\bf x}) d \sigma_d({\bf x}) \approx
\frac{1}{|X|} \sum_{{\bf x} \in X} f({\bf x})
\end{equation}
is exact for all homogeneous polynomials $f({\bf x})=f(x_0,x_1,\dots,x_d)$
of degree $t$ ($\sigma_d$ denotes the surface measure on $S^d$). A set $X$ is a
$t$-design if it is of index $k$ for every $k \leq t$.

The concept of spherical $t$-designs has been studied extensively from
various points of view, including representation theory, combinatorics, and
approximation theory. For general references see \cite{Baj:1992a} and
\cite{DGS}. The existence of spherical designs for every $t$ and $d$ and
large enough $N=|X|$ was first proved by Seymour and Zaslavsky in 1984
\cite{SeyZas:1984a}.

The first general construction of spherical designs for arbitrary $t$, $d$,
and large enough $N$ was given independently by Wagner \cite{W} and Bajnok
\cite{Baj:1992a}, who used $N \geq C(d)t^{O(d^4)}$ and $N \geq
C(d)t^{O(d^3)}$ points, respectively. This bound was later reduced to
$C(d)t^{d^2/2+d/2}$ by Korevaar and Meyers \cite{KM}. They believe that the
minimum size of a $t$-design on $S^d$ is $C(d)t^d$. Because of its
theoretical and practical importance, there is a keen interest in finding
explicit constructions for spherical $t$-designs on $N$ points with $N$
relatively small with respect to $t$ and $d$.

For explicit constructions of spherical designs it is often convenient to
use the following equivalent definition, for a proof see \cite{DGS}.

\bigskip

{\bf Lemma 2.2} {\it
A finite subset $X$ of $S^d$ is a spherical $t$-design if and only if for
every homogeneous harmonic polynomial $f$ of total degree $t$ or less}
$$\sum_{{\bf x} \in X} f({\bf x})=0.$$
\bigskip

A polynomial $f(x_0,x_1,\dots,x_d)$ is \emph{harmonic} if it satisfies
Laplace's equation $\Delta f=0$. The set of homogeneous harmonic
polynomials of degree $k$ over $S^d$ forms the vector space $Harm_{k}(S^d)$
with $$\dim Harm_{k}(S^d)=\bc{d+k}{k}-\bc{d+k-2}{k-2}.$$ In particular, for
small values of $k$ we find that $\Phi_k(S^d)$ forms a basis for
$Harm_k(S^d)$ where
\bigskip

$\Phi_1(S^d)=\{x_i | 0 \leq i \leq d\},$

\medskip

$\Phi_2(S^d)=\{x_ix_j| 0 \leq i<j \leq d\} \cup \{x_i^2-x_{i+1}^2 | 0 \leq
i \leq d-1\},$

\medskip

$\Phi_3(S^d)=\{x_ix_jx_k | 0 \leq i<j<k \leq d\} \cup \{x_i^3-3x_ix_j^2 | 0
\leq i \not =j\leq d\},$
 and

\medskip

$\Phi_4(S^d)=\{x_i^3x_j-x_ix_j^3 | 1 \leq i  <j\leq d+1 \} \\ \hspace*{1in}
\cup \{x_i^4-6x_i^2x_j^2+x_j^4 | 1 \leq i  <j\leq d+1 \} \\ \hspace*{1in}
\cup \{x_i^3x_j-3x_ix_jx_k^2 | 1 \leq i  <j\leq d+1, \\ \hspace*{2.8in} 1
\leq k \leq d+1, i \not= k, j \not= k \} \\ \hspace*{1in} \cup
\{x_ix_jx_kx_l | 1 \leq i<j<k <l \leq d+1\} $.
\bigskip

We now prove
\bigskip

{\bf Proposition 2.3}\, {\it 
For every odd positive integer $k$, odd prime $p$, and dimension $d \geq
1$, $X(d,p)$
is a spherical design of index $k$.  Furthermore, $X(d,p)$ is a spherical
3-design.}

\bigskip

{\bf Proof.} Let us consider our point set $X=X(d,p)$.  Since for any
reflection $T$ with respect to the hyperplane $x_i=0$ ($0 \leq i \leq d+1$)
we have $X^{T}=X$, the equation $$\sum_{{\bf x} \in X} f({\bf x})=0$$ will
hold for every $f=\Pi_{i=0}^{d+1} x_i^{k_i} \in Harm_{k}(S^d)$ which has at
least one variable $x_i$ with an exponent $k_i$ odd. This observation
immediately implies that our point set $X$ is a spherical design of index
$k$ for every odd integer $k$.

Furthermore, our point set $X$ also remains fixed under a reflection $U$
with respect to the hyperplane $x_i=x_j$ ($0 \leq i,j \leq d+1$), thus we
also have $$\sum_{{\bf x} \in X} f({\bf x})=0$$ for the quadratics of the
form $x_i^2-x_{i+1}^2$ ($0 \leq i \leq d-1$).  Therefore $X$ satisfies
Lemma 2.2 for $t=3$ and is thus a 3-design. $\Box$

While our construction does not, in general, give a spherical $t$-design
for $t \geq 4$, it is worth mentioning a few cases when it does. Namely, we
prove the following.
\bigskip

{\bf Proposition 2.4}\, {\it The sets $X(3,3)$, $X(2,5)$
and $X(3,5)$ form spherical 5-designs, and $X(1,7)$ is a
spherical $7$-design.}
\bigskip

{\bf Proof.}  As noted before, $X(1,7)$ is the regular octagon on the
circle $S^1$, hence a well known 7-design (see \cite{CS}).  To prove that
$X(3,3)$, $X(2,5)$ and $X(3,5)$ are spherical 5-designs, it suffices to
verify that they are of index 4; this, together with Proposition 2.3 yields
that they form 5-designs.  To prove that they are of index 4, one can
easily check that $$\sum_{{\bf x} \in X} f({\bf x})=0$$ for the polynomials
in $\Phi_4(S^d)$. $\Box$
\bigskip

{\bf Remark.}\, It is worth noting that the minimum size of a spherical
$3$-design on $S^d$ is $2d+2$, and this is achieved by the vertices of the
generalized regular octahedron.  The size of our pointset, $N(d,p)$, as
given by equation (1.2), is generally a lot larger than this; however, our
construction may prove useful for other purposes.

\section{Generalized s energy} 
\setcounter{equation}{0} 

Given a set of points $\omega_N:=\{x_{1},...,x_{N}\},\, N\geq 1$
on $S^d$ and $s>0$, we define
the $s$ energy associated with $\omega_N$ by
\begin{equation}\label{3.1}
E_d(s, \omega_N):=
\sum_{1\leq i<j\leq N}\frac{1}{||x_{i}-x_{j}||^s}.
\end{equation}
\medskip

{\it Fekete points} on $S^d$, see {\cite{KS1} and the references cited therein, 
are points that minimize the $s>0$ energy over all sets of $N$ points.
Physically, 
this represents the energy of $N$
charged particles that repel each other according to Coulomb's law.
(See also \cite{RSZ} for a discussion of Elliptic Fekete points: the case
$s=0$).
 From recent work of \cite{DG}, it is known that  points of minimal energy for $s\leq d$ are well separated in the
sense that 
\[
\lim_{N\to\infty}\frac{1}{\sigma_d(S^d)} \int_{S^d}f(y) d \sigma_d(y)=
 \lim_{N\to\infty} \frac{1}{N} \sum_{y\in \omega_N} f(y)
\]
for all continuous functions $f$ on $S^d$.
Another way to understand this is to define
for a given set of points $\omega_N$,
\[
\delta(\omega_N):={\rm inf}_{i\neq j}||x_{i,N}-x_{j,N}||;\, \, \,
\delta_N:={\rm sup}_{\omega_N\subset S^d}\delta(\omega_N).
\]
The determination of $\delta_N$ is called {\it Tammes problem} or the {\it
Spherical 
packing problem}, see \cite{CS}. It asks to maximize the smallest distance 
among $N$ points on $S^d$.
More precisely, in \cite{KS1}, the authors have shown that uniformly for 
any $s>d\geq 2$,
\begin{equation}\label{eq 3.2}
\frac{1}{C}N^{1+s/d}\leq E_d(s, \omega_N^{*})\leq C N^{1+s/d}
\end{equation}  
for any set of points $\omega_N^{*}$ minimizing $E_d(s,\omega)$ over all
sets of points
$\omega_N$.  Fixing $N$ in the above and letting $s\to\infty$ then implies
that for any $x_i\in \omega_N^{*},\, s>d$,
\begin{equation}\label{3.3}
||x_{i}-x_{j}||\geq CN^{-1/d}.
\end{equation}
If $d=2$, then it is a well known result of W. Habicht and 
B. L. van der Waerden, see \cite{HW} that  
\[
\delta_N=\left(\frac{8\pi}{\sqrt{3}}\right)^{1/2}N^{-1/2}+O(N^{-2/3}),\,
N\to \infty.
\]
Thus for $d=2$, the minimal energy $s$ problem 
reduces to the best packing problem which is the optimal choice of points
one might hope for.
For $0<s<d$, the energy integral 
\[
\int_{S^d}\int_{S^d}\frac{1}{||x-y||^s}d\sigma_d(x)d\sigma_d(y)
\]
is finite and its value is
 
$$I_{s,d}:=\frac{\Gamma((d+1)/2)\Gamma(d-s)}{\Gamma((d-s+1)/2)\Gamma(d-s/2)},
$$ where $\Gamma $ denotes the gamma function.
 Using this fact, it has been shown, see \cite{W1}
and \cite{W2}, that
for $d-2\leq s\leq d$
\begin{equation}\label{eq 3.4}
E_2(s,\omega_n^*)=\frac{1}{2}I_{s,d}N^2-R_{N,s,d}
\end{equation}
where
\begin{equation}\label{eq 3.5}
\frac{1}{C}\leq \frac{R_{N,s,d}}{N^{1+s/d}}\leq C.
\end{equation}
In fact, for $d=2$ and $0<s<2$ it is conjectured, see \cite{R}, that 
\begin{eqnarray}\label{eq 3.6}
&& \lim_{N\to\infty}\frac{R_{N,s,2}}{N^{1+s/2}}= \\
&& \nonumber
-3\left(\frac{\sqrt{3}}{8\pi}\right)^{1/2}\zeta(1/2)\sum_{n=0}^{\infty}
\left(\frac{1}{\sqrt{3n+1}}-\frac{1}{\sqrt{3n+2}}\right)\approx 0.55305,
\end{eqnarray}
where $\zeta$ denotes the Riemann zeta function defined by $ \zeta(s) =\sum _{n=1}^{\infty} \frac {1}{n^s}$.
In addition for $s=d$, it is known, see \cite{KS1} that 
\begin{equation}\label{eq 3.7}
\lim_{N\to \infty}\frac{E_d(d,\omega_N^*)}{N^2\log
N}=\frac{\Gamma((d+1)/2)}{2d\Gamma(d/2)
\Gamma(1/2)}.
\end{equation}
We first concentrate our efforts in computing 
$E_2(s, X),\, s=1, 2$ and $2.5$ for $p=3,...,211$.
The numerical data is contained
in Appendix A. 


  
\medskip

{\bf Report 3.1 ($d=2$)}
\begin{itemize}
\item[(a)] Numerical data for $s=1$: Let 
\[
e_{2,3}:=\frac{E_3(2,X)}{N^2}
\]
and write in view of $(3.4)$ and $(3.5)$
\[
\frac{R_{N,1,2}}{N^{2}}=(1/2-e_{1,2})\sim \frac{1}{N^{1/2}}
\]
and
\[
\frac{R_{N,1,2}}{N^{3/2}}=(1/2-e_{1,2})N^{1/2},
\]
where here and throughout, for non-negative sequences $a_n$ and $b_n$, $a_n \sim b_n$ if $1/C \leq a_n/b_n \leq C$. 

Then for $p=31,...,211$,
\[
0.004\leq \left|\frac{R_{N,1,2}}{N^{2}}\right|\leq 0.008
\] 
and
\[
1.21\leq \left|\frac{R_{N,1,2}}{N^{3/2}}\right|\leq 1.47.
\]
Thus it seems likely that $\frac{R_{N,1,2}}{N^{2}}$ approaches $0$ and moreover
we may take $C=1.47$ in $(3.5)$.
 
\item[(b)] Numerical data for $s=2$. Let
\[
e_{2,2}:=\frac{E_2(2,X)}{N^2\log N}
\]
and write in view of $(3.7)$, 
\[
R_{N,2,2}=\left|1/8-e_{2,2}\right|.
\]
Then for $p=31,...,211$,
\[
0.025\leq R_{N,2,2}\leq 0.205. 
\]
Thus  $e_{2,2}$  deviates from $1/8$ with a maximum error of $0.205$. 
\item[(c)] Numerical data for $s=3$. Let 
\[
e_{3,2}:=\frac{E_2(3,X)}{N^{5/2}}.
\]
Then for $p=31,...,211$,
\[
0.144\leq e_{3,2}\leq 0.259.
\]
Thus we may take $C=0.259$ in $(3.2)$.
\end{itemize}

We next computed 
$E_3(s, X),\, s=2,3$ and $3.125$ for $p=3,...,83$.
The numerical data is contained
in Appendix B. 

 
\medskip

{\bf Report 3.2 ($d=3$)}
\begin{itemize}
\item[(a)] Numerical data for $s=2$: Let 
\[
e_{1,2}:=\frac{E_2(1,X)}{N^2}
\]
and write in view of $(3.4)$ and $(3.5)$
\[
\frac{R_{N,2,3}}{N^{2}}=(1/2-e_{1,2})\sim 1.
\]
Then for $p=31,...,83$,
\[
0.01\leq |\frac{R_{N,2,3}}{N^{2}}|\leq 0.02
\] 
Thus we may take  $C=0.02$ in $(3.5)$.
\item[(b)] Numerical data for $s=3$. Let
\[
e_{3,3}:=\frac{E_3(3,X)}{N^2\log N}
\]
and write in view of $(3.7)$, 
\[
R_{N,3,3}=\left|\frac{1}{3\pi}-e_{3,3}\right|.
\]
Then for $p=31,...,83$,
\[
0.053\leq R_{N,3,3}\leq 0.091. 
\]
Thus  $e_{3,3}$  deviates from $\frac{1}{3\pi}$ with a maximum
error of $0.091$. 
\item[(c)] Numerical data for $s=3.125$. Let 
\[
e_{3,3.125}:=\frac{E_3(3.25,X)}{N^{2.04}}.
\]
Then for $p=31,...,83$,
\[
2.11\leq e_{3,3.125}\leq 3.95.
\]
Thus we may take $C=3.95$ in $(3.2)$.
\end{itemize}

Finally we considered numerically results on the spacings of our points for
$d=2$ and $d=3$.
To achieve our lower bounds, it suffices to apply $(3.2)$ together
with $(3.1)$. We did this using our results above for 
$d=2,\, s=2.5$ and $d=3,\, s=3.125$. Our data below shows 
that our points on $S^d$, $d=2$ and $d=3$ are separated numerically by a factor of order $O(1/N^{1/d+1/s})$ with a perturbation factor 
of order $O(1/N^{1/s})$ to the ideal order $O(1/N^{1/d})$.
More precisely, we have:
\medskip

{\bf Report 3.3 (spacing)}\, For any $x_i\in X$, the following hold:
\begin{itemize}
\item[(a)] For $d=2$ and $p=\{3,...,211\}$
\[
||x_{i}-x_{j}||\geq \frac{1}{(0.254)^{1/(2.5)}}\frac{1}{N^{1/2+1/(2.5)}},\,
i\neq j.
\]
\item[(b)] For $d=3$ and $p=\{3,...,81\}$
\[
||x_{i}-x_{j}||\geq
\frac{1}{(3.95)^{1/(3.125)}}\frac{1}{N^{1/3+1/(3.125)}},\, i\neq j.
\]
\end{itemize}

\section{Extension to finite fields of odd prime powers}
\setcounter{equation}{0}

In this last section, we briefly observe
that we may solve the same quadratic form (1.1) over a general 
finite field $F_q$, where $q=p^e$ is an odd prime power and in this way
distribute points on 
$S^d$ as well. One way to do this is as follows. Assume that $q=p^e$, with
$e\geq 1$. Then 
the field $F_q$ is an $e$-dimensional vector space over the field 
$F_p$. Let $\alpha _1, \dots, \alpha_e$ be a basis of $F_q$ over $F_p$.
Thus if 
$\alpha \in F_q$, then $\alpha$ can be uniquely written as 
$\alpha = a_1\alpha _1 +\dots + a_e\alpha_e$, where each $a_i \in F_p$. 
Moreover, we may assume that each $a_i$ satisfies $-(p-1)/2 \leq a_i \leq
(p-1)/2$. 

If $(x_1, \dots ,x_{d+1})$ is a solution to the quadratic form (1.1) over
$F_q$, 
then each $x_i$ is of the form $x_i=\alpha \in F_q$. Corresponding to the 
finite field element $x_i=\alpha$, we may now naturally associate the integer 
$M_i = a_1 +a_2p+\dots +a_ep^{e-1}$. It is an easy exercise to check that
indeed $ -(p^e-1)/2 \leq M _i \leq (p^e-1)/2$. 
We then map the vector $V=(M_1, \dots, M_{d+1})$ to the surface of the unit
sphere 
$S^d$ by normalizing the vector $V$. We note that when $e=1$, this reduces
to our original construction. 
In particular, for increasing values of $e$, we obtain an
increasing number $N_e$ of points scattered on the surface of the unit sphere 
$S^d$, so that as $e \rightarrow \infty$, it is clear that 
$N_e \rightarrow \infty$.

In light of the success of our initial investigation 
for the special case $e=1$,  a detailed analysis of the cases where $e > 1$  may warrant further investigation.

\vskip .2in

{\bf Acknowledgments.}
The research of (SBD) began while visiting Penn State University during
the 1999 academic year and is supported, in part, by a 2000 Georgia
Southern University 
research grant.   
(SBD) would like to thank the Mathematics Department at Penn State
University for their
kind support and hospitality. We would  like to thank Andrew Clayman for some preliminary calculations related to this project.

\bigskip

$^{1}$ Department of Mathematics, 
Gettysburg College, Gettysburg, PA 17325, Email: bbajnok@gettysburg.edu

$^{2}$ Department of Mathematics and Computer Science, Georgia Southern
University, P. O. Box 8093, Statesboro, GA 31406, Email:
damelin@gsu.cs.gasou.edu

$^{3}$,$^{4}$ Department of Mathematics, The Pennsylvania State University, 
University Park, PA  16802, Email: li@math.psu.edu, mullen@math.psu.edu

\bigskip

{\bf Appendix A},\, Table of normalized $E_2(s,X),\, s=1,2,2.5$\newline
For convenience, we adopt the following conventions: \newline
Set $e1:=e_{1,2}=\frac{E_2(1,X)}{N^2}$,
$e2:=e_{2,2}=\frac{E_2(2,X)}{N^2\log N}$,
$e3:=e_{3,2}:=\frac{E_2(3,X)}{N^{5/2}}$.
\begin{displaymath}
\begin{array}{ccccc}
p=3  & e1, e2, e3 = &   0.27736892706218 &   0.01046457424783 &
0.11980183204618 \\
p=5  & e1, e2, e3 = &    0.40195539994316 &  0.11628809346709 &
0.16900002651468 \\
p=7  & e1, e2, e3 = &   0.43294284766539 &  0.14155651348989 &
0.20783488490901 \\
p=11 & e1, e2, e3 = &   0.47263289296172 &  0.15758965144183 &
0.22872928195394 \\
p=13 & e1, e2, e3 = &   0.46734388392531 &  0.13930424776208 &
0.19737146213200 \\
p=17 & e1, e2, e3 = &   0.47831260388752 &  0.14292677178733 &
0.19559212838090 \\
p=19 & e1, e2, e3 = &   0.48862305761514 & 0.15647575327531 &
0.21230845889239 \\
p=23 & e1, e2, e3 = &   0.49680020536376 &  0.16083024644030 &
0.21114298075140 \\
p=29 & e1, e2, e3 = &   0.49963999603979 &  0.20411306813044 &
0.25438001898121 \\
p=31 & e1, e2, e3 = &   0.50558991175900 &  0.20872789826935 &
0.25835091872858 \\
p=37 & e1, e2, e3 = &   0.49679590269129 &  0.16908367785692 &
0.20014914672875 \\
p=41 & e1, e2, e3 = &   0.49967248912416 &  0.18213470230307 &
0.21067502961875 \\
p=43 & e1, e2, e3 = &   0.50082510107625 &  0.15645393314738 &
0.17997113615146 \\
p=47 & e1, e2, e3 = &   0.50395991784508 &  0.17448493068513 &
0.19648831884456 \\
p=53 & e1, e2, e3 = &   0.50302988221990 &  0.20628228220202 &
0.22447505302842 \\
p=59 & e1, e2, e3 = &   0.50770676508130 &  0.21250313388376 &
0.22610594444172 \\
p=61 & e1, e2, e3 = &   0.50495658406683 &  0.22070120695148 &
0.23184436755631 \\
p=67 & e1, e2, e3 = &   0.50654515261995 &  0.19456180151876 &
0.20028103675876 \\
p=71 & e1, e2, e3 = &   0.50799678645242 &  0.19531553592815 &
0.19798724786361 \\
p=73 & e1, e2, e3 = &   0.50138319997937 &  0.17763041043367 &
0.17807400648780 \\
p=79 & e1, e2, e3 = &   0.50667682389347 &  0.24131448167150 &
0.23767087684753 \\
p=83 & e1, e2, e3 = &   0.50632565820368 &  0.19072462271388 &
0.18532160079489 \\
p=89 & e1, e2, e3 = &   0.50420397664596 &  0.19014449465661 &
0.18065926631423 \\
p=97 & e1, e2, e3 = &   0.50484339249495 &  0.17456698846262  &
0.16193580914486 \\
p=101 & e1, e2, e3 =&    0.50561387064791 &  0.21342480676679 &
0.19574508602563 \\
p=103 & e1, e2, e3 =&    0.50800632405351 &  0.19108133195127 &
0.17476558366460 \\
p=107 & e1, e2, e3 =&    0.50744167603236 &  0.21193413555663 &
0.19173518130406 \\
p=109 & e1, e2, e3 =&    0.50530215087720 &  0.21451186060578 &
0.19252922201869 \\
p=113 & e1, e2, e3 =&    0.50508191824201 &  0.22699743097426 &
0.20164193844709 \\
p=127 & e1, e2, e3 =&    0.50685489076967 &  0.18562274384222 &
0.15976587408138 \\
p=131 & e1, e2, e3 =&    0.50784482237133 &  0.21884299920671 &
0.18664215587238 \\
p=137 & e1, e2, e3 =&    0.50588778915797 &  0.22528035319823 &
0.18918543160007 \\
p=139 & e1, e2, e3 =&    0.50868464653880 &  0.22665247011204 &
0.18992843422136 \\
p=149 & e1, e2, e3 =&    0.50654747356646 &  0.25299360855091 &
0.20721585270831 \\
p=151 & e1, e2, e3 =&    0.50839219613404 &  0.21335681640995 &
0.17440146172139 \\
p=157 & e1, e2, e3 =&    0.50550126525598 &  0.24862960416529 &
0.20046770502092 \\
p=163 & e1, e2, e3 =&    0.50708965412530 &  0.18736597295888 &
0.14964778305454 \\
p=167 & e1, e2, e3 =&    0.50726676470783 &  0.19666444095315 &
0.15591738945230 \\
p=173 & e1, e2, e3 =&    0.50594513571582 &  0.20098724096747 &
0.15735377427577  \\
p=179 & e1, e2, e3 =&    0.50824611328744 &  0.33454655513670 &
0.25964637978174 \\
p=181 & e1, e2, e3 =&    0.50687976341157 & 0.27359266111448 &
0.21125394734360 \\
p=191 & e1, e2, e3 =&    0.50841113729910 &  0.20611531648543 &
0.15679195137179 \\
p=193 & e1, e2, e3 =&    0.50641535139461 &  0.23733139042064 &
0.17966621151199 \\
p=197 & e1, e2, e3 =&    0.50665988263156 &  0.22671894926843 &
0.17054554488575 \\
p=199 & e1, e2, e3 =&    0.50721503705415 &  0.22942421796768 &
0.17230976048709 \\
p=211 & e1, e2, e3 =&    0.50723259692356 &  0.19527807424465 &
0.14400243632793
\end{array}
\end{displaymath} 

\newpage

{\bf Appendix B},\, Table of normalized $E_3(s,X),\, s=2,3,3.125$\newline
For convenience, we adopt the following conventions:

\bigskip

Set $e1:=e_{2,3}=\frac{E_3(2,X)}{N^2}$,
$e2:=e_{3,3}=\frac{E_3(3,X)}{N^2\log N}$, $e3:=e_{3,3.125}:=
\frac{E_3(3.125,X)}{N^{2.04}}$.
\begin{displaymath}
\begin{array}{ccccc}

 p=3 &  e1, e2, e3=&  0.28993055555556 & 0.07726112747611  & 0.24132173335499 \\
 p=5 &  e1, e2, e3=&  0.37650018037519 & 0.08604648073508 &  0.42113832739910 \\
 p=7 &  e1, e2, e3=&  0.42533543875492 & 0.09834709791994 & 0.60593660203289 \\
 p=11 &  e1, e2, e3=&  0.49141162848846 &0.14167115959202 &
1.1735995123910 \\
 p=13 &  e1, e2, e3=&  0.48120035766978 &0.12590020472700 &
1.1139764584208 \\
 p=17 &  e1, e2, e3 = & 0.52132037047210 &  0.19355085551622 &
2.0982770508338 \\
 p=19 & e1, e2, e3 = & 0.50870874907126 &  0.14945260382314 &
1.6099778963136 \\
 p=23 & e1, e2, e3 = & 0.50664591863877 &  0.15354023703434 &
1.8296037700500 \\
 p=29 & e1, e2, e3 = & 0.51703948372378 &  0.15825300052952 &
2.0440071303765 \\
 p=31 & e1, e2, e3 = & 0.51738403568925 &  0.15916399708017 &
2.1127707717489 \\
 p=37 & e1, e2, e3 = & 0.51884012094570 &  0.18738716021182 &
2.8402487195221 \\
 p=41 & e1, e2, e3 = & 0.52328941444401 &  0.17241718697984 &
2.5921542948948 \\
 p=43 & e1, e2, e3 = & 0.52006281540841 &  0.16742067897271 &
2.5555792369681 \\
 p=47 & e1, e2, e3 = & 0.52154686231626 &  0.19360006853089 &
3.2832846733182 \\
 p=53 & e1, e2, e3 = & 0.52302701242501 &  0.21637834813051 &
3.9535044537455 \\
 p=59 & e1, e2, e3 = & 0.52414896518911 &  0.26040121306624 &
5.3775580289229 \\
 p=67 & e1, e2, e3 = & 0.52352555470433 &  0.16356808691818 &
2.8640336098518 \\
 p=71 & e1, e2, e3 = & 0.52278118326934 &  0.16964255435064 &
3.0594908965203 \\
 p=73 & e1, e2, e3 = & 0.52466921979023 &  0.19727729124473 &
3.8343687834560 \\
 p=79 & e1, e2, e3 = & 0.52456644056252 &  0.18994184763923 &
3.7615796423715 \\
 p=83 & e1, e2, e3 = & 0.52440012994288 &  0.18298261836839 & 3.5924738028014
 \end{array}
 \end{displaymath}
\bigskip

\end{document}